\documentclass[reqno]{amsart}
\usepackage{amsfonts,amssymb,amsmath,amsthm, graphicx, color, tikz,tikz-cd,soul,phonetic,mathrsfs}
\usepackage{url}
\usepackage[normalem]{ulem}
\usepackage{enumerate}

\newtheorem{theorem}{Theorem}[section]
\newtheorem{lemma}[theorem]{Lemma}
\newtheorem{proposition}[theorem]{Proposition}
\newtheorem{corollary}[theorem]{Corollary}
\theoremstyle{definition}
\newtheorem{definition}[theorem]{Definition}
\newtheorem{rem}[theorem]{Remark}
\newtheorem{example}[theorem]{Example}
\numberwithin{equation}{section}

\newcommand{\bte}{\begin{theorem}\quad  }
\newcommand{\ete}{\end{theorem} }
\newcommand{\bpr}{\begin{proposition}\quad  }
\newcommand{\epr}{\end{proposition} }
\newcommand{\ble}{\begin{lemma}\quad }
\newcommand{\ele}{\end{lemma}}
\newcommand{\bco}{\begin{corollary}\quad }
\newcommand{\eco}{\end{corollary} }
\newcommand{\bex}{\begin{example}\quad \rm }
\newcommand{\eex}{\end{example} }
\newcommand{\bde}{\begin{defi}\quad \rm }
\newcommand{\ede}{\end{defi} }
\newcommand{\brm}{\begin{rem} \quad \rm}
\newcommand{\erm}{\end{rem} }
\newcommand{\bpf}{\begin{proof}[{\bf{Proof.\quad}}] \rm}
\newcommand{\epf}{ \end{proof}}
\newcommand{\bdm}{\begin{displaymath} }
\newcommand{\edm}{\end{displaymath} }
\newcommand{\be}{\begin{eqnarray*}}
\newcommand{\ee}{\end{eqnarray*}  }

\newcommand{\lb}{\label}

\makeatletter
\newcommand\cupop{\mathop{\operator@font \cup}\nolimits}
\makeatother
\numberwithin{equation}{section}

\begin{document}
\author{Mojtaba Sedaghatjoo}
\address{Department of Mathematics, College of Sciences, Persian Gulf University, Bushehr, Iran.}
\email{Sedaghat@pgu.ac.ir}
\author{Mohammad Roueentan}
\address{College of Engneering, Lamerd Higher Education Center, Lamerd, Iran.}
\email{rooeintan@lamerdhec.ac.ir}
\title[Right subdirectly irreducible completely 0-simple semigroups]{Right subdirectly irreducible completely 0-simple semigroups and uniform acts over rectangular bands}
\begin{abstract}
The aim of this paper is characterizing right subdirectly irreducible completely 0-simple semigroups. We prove that such semigroups are indeed groups with least nontrivial subgroups. On the other hand we prove that right irreducible completely 0-simple semigroups are groups for which nontrivial subgroups have nontrivial intersection. Ultimately, we characterize the class of uniform acts as an overclass of subdirectly irreducible acts over rectangular bands.
\end{abstract}

\maketitle
\section{INTRODUCTION AND PRELIMINARIES}
Apparently, investigating semigroups possessing least nondiagonal right congruences, termed right subdirectly irreducible semigroups, was initiated by Rankin et al. \cite{Ran}, who presented a general account on such semigroups. This class of semigroups is indeed a subclass of subdirectly irreducible semigroups on which the first investigations were pioneered by the efforts of Thierrin \cite{thier} and Schein \cite{schein}. On the other hand, representations of semigroups by transformations of sets (acts over semigroups), have been played an undeniable role in semigroup theory which a great deal of works have been devoted to this area. Thereby, in light of determining any variety of algebras of the same type by its subdirectly irreducible members (Birkhoff's theorem), the class of subdirectly irreducible acts, was a matter of interest in literatures and was followed by some authors for instance Kozhukhov in \cite{kozh1,kozh2,kozh3,kozh}. Uniform acts that are generalizations of subdirectly irreducible acts, were investigated in a specific case in \cite{feller} by Feller and Gantos as a pioneering work whilst the recent work \cite{MM}, presents a general account on this issue.

Throughout this paper, $S$ will denote a semigroup. To every semigroup $S$ we can associate the monoid $S^1$ with the identity element $1$ adjoined if necessary. Indeed, $S^1=\begin{cases} S & \text{if} ~S~ \text{has an identity element,} \\ S  \cup\{1\} & \text{otherwise}. \end{cases}$

A right $S$-act $A_S$ (or an act $A$, if there is no danger of ambiguity) is a non-empty set together with a right action $ \mu: A\times S \longrightarrow A ,\,as:=\mu(a,s)$, such that $a(st)=(as)t$ for each  $a\in A$ and  $s,t\in S$. Hereby, any semigroup $S$ can be considered a right $S$-act over itself with the natural right action, denoted by $S_S$.
An element $\theta$ of an act  $A$ is said to be a zero element if $\theta s=\theta$ for all $s\in S$. The set of zero elements of an act $A$ is denoted by $Z(A)$. Moreover, the one element act, denoted by $\Theta=\{\theta\}$, is called the zero act. Recall that an act is called simple ($\theta$-simple) if it contains no (nonzero) subacts other than itself. A semigroup $S$ is called simple ($\theta$-simple)  if $S_S$ is simple ($\theta$-simple).  An equivalence relation  $\rho$ on an $S$-act $A$ is called a right congruence or simply a congruence on $A$ if  $a~\rho ~a'$ implies $(as)~ \rho ~(a's)$ for every $a,a' \in A, s \in S$ and the class of $a$ under $\rho$ is denoted by $a_{\rho}$. The set of all congruences on $A$ is denoted by Con($A$). For an act $A$ the diagonal relation $\{(a,a)\,|\,a\in A\}$ on $A$ is a congruence on $A$ which is denoted by $\Delta_A$. Also if $B$ is a subact of $A$, then the congruence $(B\times B)\cup \Delta_A$  on $A$ is denoted by $\rho_B$ and is called the Rees congruence by the subact $B$.
For $a,b\in A, $ the monocyclic congruence on $A$ generated by the pair $a,b$ is denoted by $\rho(a,b)$ which for $x,y \in A, x \rho(a,b)y$ if and only if $x=y$ or there exist $p_1,p_2,\ldots ,p_n,q_1,q_2,\ldots ,q_n\in A, w_1,w_2,\ldots ,w_n\in S^1$ where for every $i=1,2,\ldots ,n, (p_i,q_i) \in \{(a,b), (b,a)\}$, with the following sequence of equalities:
\begin{alignat*}{3}
x =&p_1 w_1 & &q_2 w_2 =p _3
w_3 & &~\cdots ~q_n w_n =y,\\
&q_1 w_1 \,&= \,\,&p_2  w_2 & &~\cdots ~
\end{alignat*}
shall be called of length $n$. We recall that an $S$-act $A$ is called subdirectly irreducible if every set of congruences $\{\rho_i \,| \,i\in I \}$ on $A$ with $\bigcap_{_{i\in I}} {\rho_i} = \Delta_A$, contains $\Delta_A$, indeed, the set of non-diagonal congruences has the least element. It is known that each subdirectly irreducible act $A$ contains the  minimum nontrivial subact which shall be called the kernel of $A$ and shall be denoted by Ker$(A)$. Also an act $A$ is called irreducible if for any two congruences $\rho$ and $\lambda$ on $A$, $\rho \bigcap\lambda =\Delta_A$ implies that $\rho=\Delta_A$ or $\lambda=\Delta_A$, equivalently, intersection of finite non-diagonal congruences is non-diagonal. Analogous arguments and terminologies can be employed for semigroups as right acts over themselves. For a thorough account on the preliminaries, the reader is referred to \cite{How,kilp,MM}.

\section{Right subdirectly irreducible completely 0-simple semigroups}
Since a semigroup $S$ is right subdirectly irreducible semigroups if and only if $S_S$ is right subdirectly irreducible as a right $S$-act, and considering the implication  ``subdirectly irreducible act" $\Longrightarrow$ ``irreducible act" $\Longrightarrow$ ``uniform act", in this section we bring out preliminary and basic properties of uniform acts needed in the sequel.  First recall that a subact $B$ of an act $A$ is called large in $A$ (or $A$ is called an  essential extension of $B$), denoted by $B \subseteq' A$, if any $S$-homomorphism $g: A \longrightarrow C$ such that $g|_B$ is a monomorphism is itself a monomorphism. One may routinely observe that a subact $B$ of an act $A$ is large in $A$ if and only if for every (monocyclic) non-diagonal congruence $\rho \in$ Con$(A)$, $\rho_B \cap \rho \neq \Delta_A$.

\begin{definition} For a semigroup $S$, an $S$-act $A$ is called $uniform$ if every nonzero subact is large in $A$. Also a semigroup $S$ is called right (left) uniform if the right (left) $S$-act $S_S$ ($_SS$) is uniform.
\end{definition}

\begin{definition} A subset $B$ of a right $S$-act $A$ is called separated, provided that for $a\neq b \in B$ there exists $s\in S\backslash \{1\}$ such that $as\neq bs$.
\end{definition}
The next two results, stated below, are preliminary results of uniform acts needed in the next argument.
\bco \lb{co11}\cite[Corollary 3.15]{MM}
If $S$ is a right uniform semigroup and $xy=y$ for $x,y\in S$, then $x$ is a left identity or $y$ is a left zero.
\eco
\bpr \lb{pr16} \cite[Proposition 2.3]{MM} Let $S$ be a semigroup and $A$ be a nonzero uniform $S$-act. If $B,C$
are nonzero subacts of $A$, then $|B \cap C| \geq 2$.
\epr

As a result of Corollary \ref{co11}, in right uniform semigroups, idempotents are left zeros or left identities. This important result leads to clarifying the structure of right uniform completely (0-)simple semigroups. Since any completely (0-)simple semigroup is isomorphic to the Rees matrix semigroups $S= \mathcal{M}^{(0)}[G;I,\Lambda;P]$ with $\Lambda \times I$ regular sandwich matrix $P$ over the ($0$-)group $(G^0)G$, in the sequel, the term $\mathcal{M}^{(0)}[G;I,\Lambda;P]$ stands for a completely (0-)simple semigroup (see \cite{How}). Hereby, a typical idempotent $e=(i,p_{\lambda i}^{-1},\lambda)$ is a left zero or left identity. If $S$ is completely 0-simple semigroup, since $S$ contains a zero element, $e$ must be a left identity and hence $(i,p_{\lambda i}^{-1},\lambda)(j,a,\mu)=(j,a,\mu)$ which necessitates that $j=i$. Therefore $I$ is singleton or equivalently, $S$ is right 0-simple. Now the regularity of $P$ implies that all entries are nonzero and consequently, products of nonzero elements are nonzero. On the other hand, if $S$ is completely simple semigroup, the idempotent element $e=(i,p_{\lambda i}^{-1},\lambda)$ may be a left zero element. If this occurs, then for any $g\in G$, $j\in I$ and $\mu \in \Lambda$, $e=(i,p_{\lambda i}^{-1},\lambda)(j,g,\mu)=(i,p_{\lambda i}^{-1},\lambda)$ and then $\Lambda =\{\lambda\}$ and $ p_{\lambda j}=g^{-1}$ which provides that $\Lambda$ is a singleton and $G$ is trivial. Therefore $S$ is a left zero semigroup and being uniform implies that $|S|=2$.  The next proposition is an immediate result of this argument.
\bpr \lb{pr14} Let $S=\mathcal{M}^{(0)}[G;I,\Lambda;P]$ be a completely (0-)simple semigroup with more than two elements. $S$ is uniform if and only if $I$ is a singleton.
\epr

The next result is an straightforward result of the multiplication rule in Ress matrix semigroups.
\ble \lb{le1} Let $S= \mathcal{M}^0[G;I,\Lambda;P]$ and $\mathfrak{m}=(i,a, \lambda), \mathfrak{n}=(i, b,\eta),\mathfrak{p}=(i,c,\vartheta)\in S$. Then $\mathfrak{m}\mathfrak{n}=\mathfrak{m}\mathfrak{p}$ if and only if $\mathfrak{n}=\mathfrak{p}$.
\ele

In what follows we identify the structure of monocyclic right congruences over completely (0-)simple uniform semigroups in two steps, involving identifying the classes associated to generating pairs and the other classes. Indeed, let $S= \mathcal{M}^0[G;I,\Lambda;P]$ be a completely 0-simple semigroup where $I=\{i\}$ is a singleton and let $\rho$ be a right congruence on $S$. It should be mentioned that $0_{\rho}$ is a right ideal of $S$ and hence $0_{\rho}=\{0\}$ or $S$. If $\rho=\rho(0,\mathfrak{m})$ where $\mathfrak{m}\neq 0$, $0_{\rho}=S$ which implies that $\rho=\nabla_S$. So, regarding the argument before Proposition \ref{pr14} and the sequence of equalities identifying two elements for a monocyclic congruence, $\rho$ is generated by a pair of two different nonzero elements if and only if $0_{\rho}=\{0\}$.

 Now, suppose that $\rho=\rho(\mathfrak{m},\mathfrak{n})$ where $\mathfrak{m}=(i,a,\lambda)$, $\mathfrak{n}=(i,b,\mu)$ and $\mathfrak{m} \neq \mathfrak{n}$. We assign the element $ap_{\lambda i}(bp_{\mu i})^{-1}$ to $\rho$, denoted by $X_{(\mathfrak{m},\mathfrak{n})}$ or simply $X$ if there is no ambiguity. Then $\mathfrak{a}\,\rho \, \mathfrak{b}$ if and only if $\mathfrak{a}= \mathfrak{b}$ or there exist $p_1,p_2,\ldots ,p_n,q_1,q_2,\ldots ,q_n\in S, w_1,w_2,\ldots ,w_n\in S^1$ where for every $k=1,2,\ldots ,n,$ $ (p_k,q_k) \in \{(\mathfrak{m},\mathfrak{n}),(\mathfrak{n},\mathfrak{m})\}$, with
\begin{equation} \lb{f1}
\begin{alignedat}{3}
\mathfrak{a} =&p_1 w_1 & &q_2 w_2 =p _3
w_3 & &~\cdots ~q_n w_n =\mathfrak{b}.\\
&q_1 w_1 \,&= \,\,&p_2  w_2 & &~\cdots ~
\end{alignedat}
\end{equation}
 With no lose of generality we can assume that this sequence is of the minimum length. Therefore, in view of Lemma \ref{le1}, for any $2\leq k \leq n$, if $w_{k-1},w_k\neq 1$, $p_k\neq q_{k-1}$.

First we identify the class $\mathfrak{m}$ (identically $\mathfrak{n}$) by setting $\mathfrak{a}=\mathfrak{m}$. Again the minimality of the sequence length, implies that $w_k\neq 1$ for every $2\leq k \leq n$ (if $w_k=1$ for some $2 \leq k \leq n$ then $q_kw_k$ equals to $\mathfrak{m}$ or $\mathfrak{n}$ which contradicts the minimality of the length of the sequence). This property eventuate in two cases.

{\bf Case 1:\,$w_1=1$.}

As $\mathfrak{a}=\mathfrak{m}$ and $w_1=1$, $p_1=\mathfrak{m}$ and $q_1=\mathfrak{n}$. In this case two subcases may occur.

\hspace{2cm}{\bf Subcase 1: $p_2=\mathfrak{n}$.}

In this subcase $q_2=\mathfrak{m}$ and the equality $q_1w_1=p_2w_2$ yields $w_2=(i,p_{\mu i}^{-1},\mu)$. Thus $q_2w_2=\mathfrak{m} w_2=(i,ap_{\lambda i}p_{\mu i}^{-1},\mu)=(i,Xb,\mu)$. Regarding the argument after the sequence \ref{f1}, for any $2\leq k \leq n$, $p_k=\mathfrak{n}$ and $q_k=\mathfrak{m}$ and hence taking $w_3=(i,x_2,\mu)$ for some $x_2\in G$, $q_2w_2=(i,Xb,\mu)=\mathfrak{n}(i,x_2,\mu)$ for some $x_2\in G$. The last equality implies that $Xb=bp_{\mu i}x_2$, equivalently, $p_{\mu i}^{-1}b^{-1}Xb=x_2$. Thus \[q_3w_3=(i,a,\lambda)(i,x_2,\mu)=(i,ap_{\lambda i}x_2,\mu)=(i,ap_{\lambda i}p_{\mu i}^{-1}b^{-1}Xb,\mu)=(i,X^2b,\mu).\]
An inductive process reaches us to the equality $(i,X^{n-1}b,\mu)=q_nw_n=\mathfrak{b}.$ On the other hand a recursive process shows that any element of the form $(i,X^{n}b,\mu)$ for some natural $n$ is in $\mathfrak{m}_{\rho}$.

\hspace{2cm}{\bf Subcase 2: $p_2=\mathfrak{m}$.}

In this subcase, an argument parallel to subcase 1, provides that $(i,X^{-(n-1)}b,\mu) =q_nw_n=\mathfrak{b}$ and any element of the form $(i,X^{-n}b,\mu)$ for some natural $n$, is in $\mathfrak{m}_{\rho}$. Therefore, in case 1, a typical elements is in $\mathfrak{m}_{\rho}$ if and only if it has the form  $(i,X^nb,\mu)$ for some $n\in \mathbb{Z}$.

{\bf Case 2:\,$w_1\neq1$.}

For any $1\leq k \leq n$, since $w_k\neq 1$, $w_k=(i,x_k,\lambda)$ for some $x_k\in G$. Considering $p_1\in \{\mathfrak{m},\mathfrak{n}\}$, this case leads to two subcases.

\hspace{2cm}{\bf Subcase 1: $p_1=\mathfrak{n}$.}

The equality $\mathfrak{m}=p_1w_1$ implies that $x_1=p_{\mu i}^{-1}b^{-1}a$ and then \[q_1w_1=\mathfrak{m}w_1=(i,ap_{\lambda i}p_{\mu i}^{-1}b^{-1}a,\lambda)=(i,Xa,\lambda).\]
Analogous to the subcases of case 1, proceeding inductively we have $(i,X^{n}a,\lambda)=q_nw_n=\mathfrak{b}$ and any element of the form $(i,X^{n}a,\lambda)$ for some natural $n$ is in  $\mathfrak{m}_{\rho}$.

\hspace{2cm}{\bf Subcase 2: $p_1=\mathfrak{m}$.}

Analogously, we reach to the equality $(i,X^{-n}a,\lambda)=q_nw_n=\mathfrak{b}$ and a typical element of the form  $(i,X^{-n}a,\lambda)$ for some natural $n$ is in $\mathfrak{m}_{\rho}$.

By virtue of both cases we realize that \[\mathfrak{m}_{\rho}=\big\{(i,X^nb,\mu)\,|\,n\in \mathbb{Z}\big\} \cup \big\{(i,X^na,\lambda)\,|\,n\in \mathbb{Z}\big\}\]

 Now, we identify the class of an arbitrary nonzero element $\mathfrak{z}= (i,z,\vartheta)$ not involved in $\mathfrak{m}_{\rho}$. This necessitates that in the sequence \ref{f1}, $w_k\neq 1$ for any $1 \leq k \leq n$. An adaptation of `Case 2' in the above argument, yields that $\mathfrak{z}_{\rho}=\{(i,X^nz,\vartheta)\,|\,n\in \mathbb{Z}\}$.

As a result of the above argument, the next proposition is obtained.

\bpr \lb{pr15} let $S= \mathcal{M}^0[G;I,\Lambda;P]$ be a completely 0-simple semigroup where $I=\{i\}$ is a singleton. Suppose that $\rho=\rho(\mathfrak{m},\mathfrak{n})$ where $\mathfrak{m}=(i,a,\lambda)$, $\mathfrak{n}=(i,b,\mu)$ and $\mathfrak{m} \neq \mathfrak{n}$. Then $\mathfrak{a}\,\rho \, \mathfrak{b}$ if and only if $\mathfrak{a}=\mathfrak{b}$ or $\mathfrak{a},\mathfrak{b}\in \big\{(i,X^nb,\mu)\,|\,n\in \mathbb{Z}\big\} \cup \big\{(i,X^na,\lambda)\,|\,n\in \mathbb{Z}\big\}$ or  $\mathfrak{a},\mathfrak{b}\in \{(i,X^nz,\vartheta)\,|\,n\in \mathbb{Z},z\neq a,b ~{\rm or} ~\vartheta \neq \lambda, \mu\}$ where $X=ap_{\lambda i}(bp_{\mu i})^{-1}$. In the case that $X=e$ is the identity element of $G$, then $\rho=\Delta_S \cup \big\{(\mathfrak{m},\mathfrak{n}), (\mathfrak{n},\mathfrak{m})\big\}$.
\epr

Recall that a group is called cocyclic if it has a cogenerator element or equivalently it has the least nontrivial subgroup. A cocyclic group with a zero element externally adjoined shall be called a cocyclic 0-group. In the next theorem right subdirectly irreducible completely 0-simple semigroups are characterized.

\bte \lb{th1}  Let $S$ be a completely 0-simple semigroup with more than two elements. $S$ is right subdirectly irreducible if and only if it is a cocyclic 0-group.
\ete
\bpf As for any group there exists a one to one and order preserving correspondence between  right congruences and subgroups, \cite[Example 2.8]{MM}, we just need to prove the necessity part.

Suppose that $S= \mathcal{M}^0[G;I,\Lambda;P]$ is a Rees matrix semigroup over the $0$-group $G^0$ with the regular sandwich matrix $P$ on $\Lambda \times I$. Since $S$ is right uniform, $I=\{i\}$ is a singleton. First we prove that $|\Lambda|< 3$. Suppose, contrary to our claim, that  $\lambda$, $\mu$ and $\vartheta$ are three distinct elements in $\Lambda$. Set $\mathfrak{m}=(i,e, \lambda), \mathfrak{n}=(i,p_{\lambda i} p_{\mu i}^{-1}, \mu),\mathfrak{m'}=(i,e, \lambda),\mathfrak{n'}=(i,p_{\lambda i} p_{\vartheta i}^{-1}, \vartheta)$. It is clearly seen that  $X_{(\mathfrak{m},\mathfrak{n})}=X_{(\mathfrak{m'},\mathfrak{n'})}$ is the identity element of $G$. Since $S$ is right subdirectly irreducible, we have $\rho (\mathfrak{m},\mathfrak{n})=\big\{(\mathfrak{m},\mathfrak{n}), (\mathfrak{n},\mathfrak{m})\big\}\cup \Delta_S=\rho (\mathfrak{m'},\mathfrak{n'})=\big\{(\mathfrak{m'},\mathfrak{n'}), (\mathfrak{n'},\mathfrak{m'})\big\}\cup \Delta_S$ which implies that $\big\{(\mathfrak{m},\mathfrak{n}), (\mathfrak{n},\mathfrak{m})\big\}=\big\{(\mathfrak{m'},\mathfrak{n'}), (\mathfrak{n'},\mathfrak{m'})\big\}$, a contradiction.
 By way of contradiction suppose that $\lambda$ and $\mu$ are two distinct elements in $\Lambda$. Setting $\mathfrak{m}=(i,e, \lambda), \mathfrak{n}=(i,p_{\lambda i} p_{\mu i}^{-1}, \mu),\mathfrak{m'}=(i,p_{\mu i}p_{\lambda i}^{-1}, \lambda),\mathfrak{n'}=(i,e, \mu)$ and regard to $X_{(\mathfrak{m},\mathfrak{n})}=X_{(\mathfrak{m'},\mathfrak{n'})}=e$  we again reach to $\big\{(\mathfrak{m},\mathfrak{n}), (\mathfrak{n},\mathfrak{m})\big\}=\big\{(\mathfrak{m'},\mathfrak{n'}), (\mathfrak{n'},\mathfrak{m'})\big\}$ which yields  $\mathfrak{m}=\mathfrak{m'}$ and $\mathfrak{n}=\mathfrak{n'}$. Thereby, $p_{\mu i}p_{\lambda i}^{-1}=e$ or equivalently $p_{\lambda i}=p_{\mu i}$. In this case, take $a\in G$ and set $\mathfrak{m}=(i,e, \lambda), \mathfrak{n}=(i,e, \mu),\mathfrak{m'}=(i,a, \lambda),\mathfrak{n'}=(i,a, \mu) $. Analogously, $X_{(\mathfrak{m},\mathfrak{n})}=X_{(\mathfrak{m'},\mathfrak{n'})}=e$ implies that  $\big\{(\mathfrak{m},\mathfrak{n}), (\mathfrak{n},\mathfrak{m})\big\}=\big\{(\mathfrak{m'},\mathfrak{n'}), (\mathfrak{n'},\mathfrak{m'})\big\}$ which necessitates that $a=e$. Thus $G$ is the trivial group and then $|S|=2$, contradicting our assumption. Therefore $\Lambda=\{\lambda\}$ is a singleton and it can be routinely checked that the assignment $g\to (i,p_{\lambda i}^{-1}g,\lambda)$ is an isomorphism from $G$ to $S\backslash 0$ and hence $S$ is right subdirectly irreducible if and only if $G$ is right subdirectly irreducible. Now, thanks to the structure of right congruences over groups, $G$ is right subdirectly irreducible if and only if it is cocyclic.

\epf

\brm \lb{rm1} Since the least non-diagonal right congruence is indeed a congruence, Theorem \ref{th1}, presents conditions under which completely 0-simple semigroups are subdirectly irreducible.
\erm

Replacing ``subdirectly irreducible" with ``irreducible" in the proof of Theorem \ref{th1} yields the next theorem.
\bte \lb{th4} Let $S$ be a completely 0-simple semigroup with more than two elements. $S$ is right irreducible if and only if $S$ is a group in which any two nontrivial subgroups have nontrivial intersection.
\ete

\section{subdirectly irreducible and uniform acts over rectangular bands}
This section is allocated to the structure of subdirectly irreducible and uniform acts over rectangular bands. Although, we have characterized, monocyclic right congruences on completely 0-simple semigroups in the case that $I$ is a singleton, it seems the problem of characterizing subdirectly irreducible acts over completely 0-simple semigroups needs more evidences to be dealt with. However, in \cite{kozh}, there is a characterization of subdirectly irreducible class of acts over rectangular bands, in what follows we clarify the structure of such acts explicitly. Thereafter, we characterize an overclass namely uniform acts over rectangular bands. The next proposition is an immediate result of Proposition \ref{pr14}.
\begin{proposition}\lb{pr1}
A rectangular band with more than 2 elements is right uniform if and only if it is a right zero semigroup.
\end{proposition}
The next lemma is an straightforward consequence of the basic properties of a rectangular band action.
\ble \lb{le2} Let $S=I\times \Lambda$ be a rectangular band for nonempty sets $I$, $\Lambda$ and $A$ be a subdirectly irreducible act over $S$. If $x=a(i,\lambda)\neq y=a(i,\mu)$ for some $a\in A,i\in I,\lambda, \mu \in \Lambda$, then $\rho(x,y)=\{(x,y),(y,x)\}\cup \Delta_A$ is the least non-diagonal congruence on $A$.
\ele
\bpr \lb{pr7} Let $S=I\times \Lambda$ be a rectangular band for nonempty sets $I$, $\Lambda$ and $A$ be an act over $S$ for which $Z(A)=\emptyset$. $A$ is subdirectly irreducible if and only if $A$ is a simple act of order 2.
\epr
\bpf
We just need to prove the necessity part. Let $A$ be a subdirectly irreducible act on $S$, $a\in A$ and $(i,\lambda)\in S$. Since $a(i,\lambda)$ is not a zero element, there exists $\mu \in \Lambda$ such that $a(i,\lambda)\neq a(i,\mu)$. Set $a(i,\lambda)=x, a(i,\mu)=y, \Lambda_1=\{\vartheta\in \Lambda\,|\,a(i,\vartheta)=x\},\Lambda_2=\{\vartheta\in \Lambda\,|\,a(i,\vartheta)=y\}$. Applying Lemma \ref{le2}, $\rho(x,y)=\{(x,y),(y,x)\}\cup \Delta_A$ is the least non-diagonal congruence on $A$. Now, let $b\in A, j\in I,\vartheta\in \Lambda$ and $b(j,\vartheta)\neq b(j, \gamma)$ for some $\gamma \in \Lambda$. Since Lemma \ref{le2} implies that $(b(j,\vartheta),b(j, \gamma)),(b(j,\gamma), b(j, \vartheta))\cup \Delta_A=\rho(x,y)$ is the least non-diagonal congruence on $A$, $b(j,\vartheta)\in\{x,y\}$. If $b(j,\vartheta)=x=a(i,\lambda)$, then $a(i,\vartheta)=a(i,\lambda)(i,\vartheta)=b(j,\vartheta)(i,\vartheta)=b(j,\vartheta)=x$ and hence $\vartheta\in \Lambda_1$. If $b(j,\vartheta)=y$, analogously, we conclude that $\vartheta \in \Lambda_2$. Therefore, for any $b\in A,j\in I,\vartheta \in \Lambda$, $b(j,\vartheta)=x \text{~or} ~y$ and then $\vartheta \in \Lambda_1\text{~or~} \Lambda_2$ respectively. Since $\vartheta$ was an arbitrary element in $\Lambda$, $\{\Lambda_1, \Lambda_2\}$ is a partition of $\Lambda$. Indeed, we have proved implicitly that for any $b\in A, b(I\times \Lambda_1)=x$ and $b(I\times \Lambda_2)=y$ and hence $AS=\{x,y\}$. Therefore, for any two distinct elements $\alpha, \beta\in A$, $\alpha s=\beta s$ for any $s\in S$ and hence $\{(\alpha,\beta),(\beta,\alpha)\cup \Delta_A$ is a non-diagonal congruence on $A$ which yields that $\{\alpha,\beta\}=\{x,y\}$. Therefore $A=\{x,y\}$, $xS=yS=A$ which yield that $A$ is a simple act of order 2.

\epf

Note that the above proposition, implicitly states that a subdirectly irreducible act $A$ with no zero element over a rectangular band $S=I\times \Lambda$, is of the form $\{x,y\}$ with $A(I\times \Lambda_1)=x$ and  $A(I\times \Lambda_2)=y$ where $\{\Lambda_1, \Lambda_2\}$ is a partition of $\Lambda$. In the next proposition we prove that subdirectly irreducible acts with one zero element over rectangular bands have a similar structure.
\bpr \lb{pr8} Let $S=I\times \Lambda$ be a rectangular band for nonempty sets $I$, $\Lambda$ and $A$ be an act over $S$ for which $Z(A)=\{\theta\}$. $A$ is subdirectly irreducible if and only if $A$ is a 0-simple act of order 2 or      Ker$(A)=\{x,y\}$ for nonzero elements  $x,y\in A$ such that for any element $a\in A,~a(I\times \Lambda_1)\subseteq \{x, \theta\},a(I\times \Lambda_2)\subseteq \{y, \theta\}$ (if $a\neq x,y$ we have the equality) where $\{\Lambda_1, \Lambda_2\}$ is a partition of $\Lambda$ and any two elements subset of $S$ different from $\{x,y\}$ is separated.
\epr
\bpf {\bf Necessity:} If $A$ contain a unique zero element, say $\theta$, two cases may occur:

{\bf Case 1:} $AS=\{\theta \}$.
In this case for any two distinct elements $a$ and $b$ in $A$, $\{(a,b),(b,a)\}\cup \Delta_A$ is a non-diagonal  congruence on $A$ and hence $A$ is a 0-simple act of order 2.

{\bf Case 2:} For some $a\in A$, $aS\neq \Theta$. Assume that $x=a(i,\lambda)\neq a(i,\mu)=y$ for some $i\in I,\lambda,\mu\in \Lambda$ and set $\Lambda_1=\{\vartheta\in \Lambda\,| \,a(i,\vartheta)=x\}, \Lambda_2 =\{\vartheta\in \Lambda\,|\,a(i,\vartheta)=y\}$. Therefore $\rho(x,y)=\{(x,y),(y,x)\}\cup \Delta_A$ is the least non-diagonal congruence on $A$. If for some nonzero element $b$, $bS=\Theta$ then, $\rho(b,\theta)=\{(b,\theta),(\theta,b)\}\cup \Delta_A$ is a non-diagonal congruence for which $\rho(x,y)\nsubseteq \rho(b,\theta)$, a contradiction. So, for any nonzero element $b$ in $A$, $bS\neq \Theta$. Now, suppose that for $b \in A,j\in I,\nu \in \Lambda$, $b (j, \nu)\neq \theta$, therefore $b(j,\vartheta)\neq \theta$ for any element $\vartheta \in \Lambda$ and  analogous to the proof of Proposition \ref{pr7}, we can prove that $\{\Lambda_1, \Lambda_2\}$ is a partition of $\Lambda$ and  $b(j\times \Lambda_1)=x$ and $b(j\times \Lambda_2)=y$. Therefore, substituting $b$ with $x$ and $y$ we get $\{x,y\}(I\times \Lambda_1)=\{x\},\{x,y\}(I\times \Lambda_2)=\{y\}$ ($xS=yS=\{x,y\}$). Therefore, Ker$(A)=\{x,y\}$  and for any nonzero element $a\in A,  a(I\times \Lambda_1)\subseteq \{x, \theta\},a(I\times \Lambda_2)\subseteq \{y, \theta\}$. On the other hand, for $a\neq x,y$ if $as\neq \theta$ for all $s\in S$, then $\rho(a,x)=\{(a,x),(x,a)\}\cup \Delta_A$ which leads to $\rho(x,y)\nsubseteq \rho(a,x)$, a contradiction. Thus $a(i,\lambda)=\theta$ for some $i\in I, \lambda \in \Lambda$ and hence $a(i,\vartheta)=\theta$ for any $\vartheta \in \Lambda$. So  $a(I\times \Lambda_1)=\{x,\theta\},a(I\times \Lambda_2)=\{y,\theta\}$. Now, on the contrary suppose that $\alpha ,\beta$ are a pair of elements different from $x,y$ such that $\alpha s=\beta s $ for any $s\in S$. Thus $\rho(\alpha,\beta)=\{(\alpha,\beta),(\beta,\alpha)\}\cup \Delta_A$ is a non-diagonal congruence for which $\rho(x,y)\nsubseteq \rho(\alpha,\beta)$, a contradiction.

{\bf Sufficiency:} Let $A$ be an act with more than two elements satisfying the second condition. Let $a$ be a nonzero element in $A$ such that $aS=\Theta$. Then $aS^1=\{a,\theta\}$ is a non-trivial subact of $A$ for which Ker$(A)\nsubseteq aS^1$, a contradiction. Then $aS$ is a non-trivial subact of $A$ containing $x$ and $y$. Let $a,b$ be two distinct elements such that $\{a,b\}\neq \{x,y\}$. Therefore, there exists $s\in S$ such that $as\neq bs$. So, our assumption implies that $as=\theta$ or $bs=\theta$. With no lose of generality suppose that $bs=\theta$. Thus $as\neq \theta$ and hence $as=x$ or $as=y$. Without lose of generality suppose that $as=x$. Then, $(x,\theta)\in \rho(a,b)$. Since $x$ is a nonzero element, $xt=y$ for some $t\in S$ and hence $(y,\theta)\in \rho (a,b)$ which imply that $(x,y)\in \rho(a,b)$. Therefore $\rho(x,y)$ is the least non-diagonal congruence on $A$.
\epf

Note that the above proposition states that for any subdirectly irreducible act $A$ over a rectangular band $S$ with more that two elements and $Z(A)=\{\theta\}$, there exists a partition of $S$ into two left ideals $I$, $J$ for which $aI=\{x,\theta \}$ and $aJ=\{y,\theta\}$ for any nonzero element $a\in A$, where Ker$A=\{x,y\}$ for nonzero elements $x,y\in S$. Therefore $aS=$Ker$A$ or $aS=$Ker$A\cup \Theta$ and hence $AS=$Ker$A\cup \Theta$.

Thanks to \cite[Theorem 1]{kozh} an act $A$ over a semigroup $S$ with $Z(A)=\{\theta_1, \theta_2\}$ is subdirectly irreducible if and only if for any two distinct elements $a,b\in A$, there exists $s\in S$ such that $\{as,bs\}=\{\theta_1, \theta_2\}$.  The next corollary gives extra property on such acts over rectangular bands.
\bpr \lb{pr9} Let $S=I\times \Lambda$ be a rectangular band for nonempty sets $I$, $\Lambda$ and $A$ be a subdirectly irreducible act over $S$ for which $Z(A)=\{\theta_1, \theta_2\}$. Then for any nonzero element $a\in A$, $aS=\{\theta_1, \theta_2\}$ and $|A|\leq |2^I|$.
\epr
\bpf  Clearly $\rho(\theta_1,\theta_2)=\{(\theta_1,\theta_2),(\theta_2, \theta_1)\}\cup \Delta_A$ is the least non-diagonal congruence on $A$ and Ker$(A)= \{\theta_1 ,\theta_2\}$. If for some nonzero element $a\in A$, $aS \nsubseteq \{\theta_1,\theta_2\}$, then for some $i\in I,\lambda,\mu\in \Lambda$, $x=a(i,\lambda)\neq a(i,\mu)=y$ where $x$ and $y$ are nonzero elements and $\rho(x,y)=\{(x,y),(y,x)\}\cup \Delta_A$ is a non-diagonal congruence on $A$ such that $\rho(\theta_1,\theta_2)\nsubseteq \rho(x,y)$, a contradiction. Then for any nonzero element $a\in A$, $aS\subseteq \{\theta_1,\theta_2\}$. On the other hand, if if for a nonzero element $a\in A$, $aS=\{\theta_i\}$ for some $1\leq i \leq 2$, then $\{a,\theta_i\}$ is a nontrivial subact of $A$ which contradict Ker$A=\{\theta_1,\theta_2\}$. Therefore $aS=\{\theta_1,\theta_2\}$. Now for any $a\in A$ assign the set $J_a=\{s\in S\,|\, as=\theta_1\}$. Clearly $J_{\theta_1}=S, J_{\theta_2}=\emptyset$ and $J_a$ is a right ideal of $S$ for any nonzero element $a\in A$. On the other hand, taking the empty set as a right ideal, there is a one to one order preserving correspondence between right ideals of $S$ and  subsets of $I$ ($J \longmapsto I_J=\{i\in I\,|\, (i,\lambda)\in J \text{~for some~} \lambda\in \Lambda \}$). Since $A$ is subdirectly irreducible, for any two distinct elements $a,b$, $J_a\neq J_b$ and hence the assignment $a \longmapsto I_{J_a}$ is one to one and we are done.
\epf
Note that in \cite[Corollary 3.21]{MM}, it is proved that for any subdirectly irreducible act $A$ with two zero elements over a semigroup $S$, $|A|\leq 2^{|S|+1}$. Therefore, Proposition \ref{pr9} sharpens this bound. The next theorem is an immediate consequence of propositions \ref{pr7}, \ref{pr8} and \cite[Theorem 1]{kozh}.
\bte
Let $S=I\times \Lambda$ be a rectangular band for nonempty sets $I$, $\Lambda$ and $A$ be an act over $S$ with $|A|>2$. $A$ is subdirectly irreducible if and only if it satisfies one of the following conditions:
\begin{enumerate}[{\rm i)}]
\item $Z(A)=\{\theta\}$ and Ker$(A)=\{x,y\}$ for nonzero elements  $x,y\in A$ such that for any element $a\in A,~a(I\times \Lambda_1)\subseteq \{x, \theta\},a(I\times \Lambda_2)\subseteq \{y, \theta\}$ (if $a\neq x,y$ we have the equality) where $\{\Lambda_1, \Lambda_2\}$ is a partition of $\Lambda$ and for any pair of elements $a,b$ different from the pair $x,y$, there exists $s\in S$ such that $as\neq bs$.
\item $Z(A)=\{\theta_1,\theta_2\}$ and for any two distinct elements $a,b\in A$, there exists $s\in S$ such that $\{as,bs\}=\{\theta_1, \theta_2\}$.

\end{enumerate}
\ete
In what follows we are going to characterize uniform acts over rectangular bands. The next result states that uniform acts are separated or have a kernel.
\ble Let $A$ be a uniform act over a semigroup $S$. If $A$ is not separated then $A$ has the kernel.
\bpf Suppose that for distinct elements $a,b\in A$, $as=bs$ for any $s\in S$. Then $\rho(a,b)=\{(a,b),(b,a)\}\cup \Delta_A$ is a non-diagonal congruence on $A$ which uniformness implies that $a,b\in B$ for any non-trivial subact $B$ of $A$.
\epf
\ele

\bpr \lb{pr6} Let $S=I\times \Lambda$ be a rectangular band for nonempty sets $I$, $\Lambda$ and $A$ be an act over $S$ for which $Z(A)=\emptyset$. $A$ is uniform if and only if Ker$A=aS$ for any $a\in A$ and for any pair of distinct elements $x$ and $y$ not simultaneously in Ker$A$, there exists $s\in S$ such that $xs\neq ys$.
\epr
\bpf
{\bf Necessity:} Let $A$ be a uniform act over $S$ and $a_0\in A$. For $i\in I,\lambda, \mu \in \Lambda$ assume that $x_0=a_0(i,\lambda)\neq a_0(i,\mu)=y_0$. Therefore, $x_0S=y_0S$ and $\rho(x_0,y_0)=\{(x_0,y_0),(y_0,x_0)\}\cup \Delta_A$ is a non-diagonal congruence on $A$ which necessitates that $x_0,y_0\in B$ for any nonzero subact $B$ of $A$ and hence Ker$A=x_0S=y_0S$. Substituting $a_0$ with an arbitrary element $a\in A$ in the above argument we observe that for any $s\in S, as \in$Ker$A$ and hence $aS=$Ker$A$. Now take distinct elements  $x$ and $y$ in $A$ for which $\{x,y\}\nsubseteq$Ker$A$. On the contrary suppose that $xs=ys$ for any $s\in S$. Thus $\rho(x,y)=\{(x,y),(y,x)\}\cup \Delta_A$ is a non-diagonal congruence on $A$ which contradicts the largeness of Ker$A$.

{\bf Sufficiency:} Regarding \cite[Lemma 2.15]{MM}, we just need to prove that Ker$A$ is large in $A$. Take the monocyclic congruence $\rho(a,b)$ on $A$ for two distinct elements $a$ and $b$ in $A$. If $a,b\in$Ker$A$, then $\rho(a,b)\cap \rho_{\text{Ker}A}\neq \Delta_A$
and we are done. Otherwise, our assumption implies the existence of $s\in S$ for which $as$ and $bs$ are two distinct elements in Ker$A$ and $(as,bs)\in \rho(a,b)\cap \rho_{\text{Ker}A}$.
\epf
\bpr \lb{pr10} Let $S=I\times \Lambda$ be a rectangular band for nonempty sets $I$, $\Lambda$ and $A$ be an act over $S$ for which $Z(A)=\{\theta\}$. $A$ is uniform if and only if $A=\{a,\theta \}$ where $aS=\Theta$ or for any nonzero element $a\in A,aS=$Ker$A$ or $aS=$Ker$A\cup \Theta$ and for any pair of distinct elements $x$ and $y$ which are not simultaneously in Ker$A$, there exists $s\in S$ such that $xs\neq ys$.
\epr
\bpf {\bf Necessity:} If $A$ contains a unique zero element, say $\theta$, two cases may occur:

{\bf Case 1:} $AS=\{\theta \}$.
In this case for any two distinct nonzero elements $a$ and $b$ in $A$, $\{a,\theta\},\{b,\theta\}$ are two nonzero subacts of $A$ which according to Proposition \ref{pr16}, $a=b$ and hence $A=\{a,\theta\}$ is  a 0-simple act of order 2 and $aS^1=A=$Ker$A$.

{\bf Case 2:} For some $a\in A$, $aS\neq \Theta$. Assume that $x=a(i,\lambda)\neq a(i,\mu)=y$ for some $i\in I,\lambda,\mu\in \Lambda$. Therefore $\rho(x,y)=\{(x,y),(y,x)\}\cup \Delta_A$ is a non-diagonal congruence on $A$ which necessitates that $x,y\in B$ for any nonzero subact $B$ of $A$ and hence Ker$A=xS=yS$. If for some nonzero element $b$, $bS=\Theta$ then, $\{b,\theta\}$ is a nonzero subact of $A$ which contradicts  Ker$A$ being minimum. So, for any nonzero element $b$ in $A$, $bS\neq \Theta$ and hence $bs\neq \theta $ for some $s\in S$. Now, replacing an arbitrary element $b$ with $a$ in the above argument, we conclude that $bs\neq \theta$ necessitates that $bs\in$Ker$A$. Therefore, for any nonzero element $b\in A$, $bS=$Ker$A$ or $bS=$Ker$A\cup \Theta$. The rest of the necessity part follows analogous to that of Proposition \ref{pr6}.

{\bf Sufficiency:} Let $A$ be a uniform act on $S$ with a unique zero element, say $\theta$.

Regarding \cite[Lemma 2.15]{MM}, we just need to prove that Ker$A$ is large in $A$. Take two distinct elements $a$ and $b$ in $A$ and the monocyclic congruence $\rho(a,b)$. If $a,b\in$Ker$A$, then $\rho(a,b)\cap \rho_{\text{Ker}A}\neq \Delta_A$ and we are done. Otherwise, our assumption implies the existence of $s\in S$ for which $as$ and $bs$ are two distinct elements. Two cases may occur.

{\bf Case 1:} Both $as$ and $bs$ are nonzero and consequently in Ker$A$. So $(as,bs)\in \rho(a,b)\cap \rho_{\text{Ker}A}$ as desired.

{\bf Case 2:} $as$ or $bs$ is the zero element. Without loss of generality suppose that $bs=\theta$. Since $as$ is nonzero, $asS=$Ker$A$ and then Ker$A\subseteq \theta_{\rho}$ which shows that $\rho(a,b)\cap \rho_{\text{Ker}A}\neq \Delta_A$ and we are done.

\epf

Since uniform acts with two zero elements are subdirectly irreducible, the next theorem is eventuated from the previous two propositions.

\bte \lb{th3}  Let $S=I\times \Lambda$ be a rectangular band for nonempty sets $I$, $\Lambda$ and $A$ be an act over $S$ with more than two elements. $A$ is uniform if and only if it satisfies one of the following conditions:
\begin{enumerate}[{\rm i)}]
\item if Ker$A=aS$ for any $a\in A$ and for any pair of distinct elements $x$ and $y$ not simultaneously in Ker$A$, there exists $s\in S$ such that $xs\neq ys$,
\item $Z(A)=\{\theta\}$ and for any nonzero element $a\in A,aS=$Ker$A$ or $aS=$Ker$A\cup \Theta$ and for any pair of distinct elements $x$ and $y$ which are not simultaneously in Ker$A$, there exists $s\in S$ such that $xs\neq ys$,
\item $AS=Z(A)=\{\theta_1, \theta_2\}$ and $A$ is separated.
\end{enumerate}

\ete

\end{document}